\documentclass[reqno,12pt]{amsart}
\usepackage{amssymb, latexsym, euscript}
\newtheorem{theorem}{Theorem}[section]
\newtheorem{proposition}[theorem]{Proposition}

\newtheorem{lemma}[theorem]{Lemma}

\newtheorem{corollary}[theorem]{Corollary}

\theoremstyle{definition}

\newtheorem{remark}[theorem]{Remark}

\numberwithin{equation}{section}

\begin{document}
\title[On self-adjoint operators in Krein spaces]{On self-adjoint operators in Krein spaces constructed
by Clifford algebra ${\mathcal C}l_2$}
\author[S.~Kuzhel]{Sergii~Kuzhel}
\author[O.~Patsiuk]{Oleksii~Patsiuk}

\address{Department of Applied Mathematics \\ AGH University of Science and Technology \\
30-059 Krakow, Poland} \email{kuzhel@mat.agh.edu.pl}

\address{Institute of Mathematics of the National
Academy of Sciences of Ukraine \\
3 Tere\-shchenkivska Street, 01601, Kiev-4 \\
Ukraine} \email{patsuk86@inbox.ru}


\maketitle
\begin{abstract}
Let $J$ and $R$ be anti-commuting fundamental symmetries in a Hilbert space $\mathfrak{H}$. The operators
$J$ and $R$ can be interpreted as basis (generating) elements of the complex Clifford
algebra ${\mathcal C}l_2(J,R):=\mbox{span}\{I, J, R, iJR\}$.
An arbitrary non-trivial fundamental symmetry from ${\mathcal C}l_2(J,R)$ is determined by the formula
$J_{\vec{\alpha}}=\alpha_{1}J+\alpha_{2}R+\alpha_{3}iJR$, where ${\vec{\alpha}}\in\mathbb{S}^2$.

Let $S$ be a symmetric operator that commutes with ${\mathcal C}l_2(J,R)$.
The purpose of this paper is to study the sets $\Sigma_{{J_{\vec{\alpha}}}}$ ($\forall{\vec{\alpha}}\in\mathbb{S}^2$) of
self-adjoint extensions of $S$ in  Krein spaces generated by fundamental symmetries ${{J_{\vec{\alpha}}}}$
(${{J_{\vec{\alpha}}}}$-self-adjoint extensions). We show that the sets $\Sigma_{{J_{\vec{\alpha}}}}$ and $\Sigma_{{J_{\vec{\beta}}}}$ are unitarily equivalent for different ${\vec{\alpha}}, {\vec{\beta}}\in\mathbb{S}^2$ and describe in detail
the structure of operators $A\in\Sigma_{{J_{\vec{\alpha}}}}$ with empty resolvent set.
\end{abstract}

\section{Introduction}
Let $\mathfrak{H}$ be a Hilbert space with inner product
$(\cdot,\cdot)$ and with non-trivial fundamental symmetry $J$ (i.e., $J=J^*$,
$J^2=I$, and $J\not={\pm{I}}$).

The space $\mathfrak{H}$ endowed with the indefinite inner product
(indefinite metric) \ $[\cdot,\cdot]_J:=(J{\cdot}, \cdot)$ is
called  a Krein space  $(\mathfrak{H}, [\cdot,\cdot]_J)$.

An operator $A$ acting in $\mathfrak{H}$ is called
${J}$-self-adjoint if $A$ is self-adjoint with respect to the indefinite metric
$[\cdot,\cdot]_J,$ i.e., if  $A^*{J}={J}A$.

In contrast to self-adjoint operators in Hilbert spaces (which
necessarily have a purely real spectrum), a $J$-self-adjoint operator $A$,
in general, has spectrum which is only symmetric
with respect to the real axis.
In particular, the situation where $\sigma(A)=\mathbb{C}$ (i.e., $A$ has empty resolvent set $\rho(A)=\emptyset$)
is also possible and it may indicate on a special structure of $A$.
To illustrate this point we consider a simple symmetric\footnote{with respect to the initial inner product $(\cdot,\cdot)$}
operator $S$ with deficiency indices $<2,2>$ which commutes with $J$:
$$
SJ=JS.
$$

It was recently shown \cite[Theorem 4.3]{KT} that the existence at least one $J$-self-adjoint extension $A$ of $S$ with empty resolvent set  \emph{is equivalent} to the existence of an additional fundamental symmetry $R$ in $\mathfrak{H}$ such that
\begin{equation}\label{es545}
SR=RS, \qquad JR=-RJ.
\end{equation}

The fundamental symmetries $J$ and $R$ can be interpreted as basis (generating) elements of the complex Clifford
algebra
${\mathcal C}l_2(J,R):=\mbox{span}\{I, J, R, iJR\}$ \cite{GK}. Hence, the existence of $J$-self-adjoint extensions of $S$ with empty resolvent set \emph{is equivalent} to the commutation of $S$ with an arbitrary element of the Clifford algebra ${\mathcal C}l_2(J,R)$.

In the present paper we investigate nonself-adjoint extensions of a densely defined symmetric operator $S$ assuming that
$S$ commutes with elements of ${\mathcal C}l_2(J,R)$. Precisely, we show that
an arbitrary non-trivial fundamental symmetry ${{J_{\vec{\alpha}}}}$ constructed in terms of  ${\mathcal C}l_2(J,R)$
is uniquely determined by the choice of vector ${\vec{\alpha}}$ from the unit sphere $\mathbb{S}^2$ in $\mathbb{R}^3$
(Lemma \ref{ll1}) and we study various collections $\Sigma_{{J_{\vec{\alpha}}}}$ of
${{J_{\vec{\alpha}}}}$-self-adjoint extensions of $S$.
Such a `flexibility' of fundamental symmetries is inspirited by the application to $\mathcal{PT}$-symmetric
quantum mechanics \cite{BE}, where $\mathcal{PT}$-symmetric  Hamiltonians are not necessarily can be realized
as $\mathcal{P}$-self-adjoint operators  \cite{AH, MO}. Moreover, for certain models \cite{GK},
the corresponding  $\mathcal{PT}$-symmetric operator realizations can be interpreted as
${{J_{\vec{\alpha}}}}$-self-adjoint operators when $\vec{\alpha}$ runs $\mathbb{S}^2$.

We show that the sets $\Sigma_{{J_{\vec{\alpha}}}}$ and $\Sigma_{{J_{\vec{\beta}}}}$ are unitarily equivalent for different
${\vec{\alpha}}, {\vec{\beta}}\in\mathbb{S}^2$ (Theorem \ref{es12a}) and describe properties of
$A\in\Sigma_{{J_{\vec{\alpha}}}}$ in terms of boundary triplets (subsections 2.4, 2.5).

Denote by $\Xi_{{\vec{\alpha}}}$ the collection of all operators $A\in\Sigma_{J_{\vec{\alpha}}}$ with empty resolvent set.
It follows from our results that, as a rule, an operator $A\in\Xi_{{\vec{\alpha}}}$ is $J_{\vec{\beta}}$-self-adjoint
(i.e., $A\in\Sigma_{{J_{\vec{\alpha}}}}\cap\Sigma_{{J_{\vec{\beta}}}}$) for a special choice of  $\vec{\beta}\in\mathbb{S}^2$  which
depends on $A$. In this way, for the case of symmetric operators $S$ with deficiency indices $<2,2>$, the complete
description of $\Xi_{{\vec{\alpha}}}$ is obtained as the union of
operators $A\in\Sigma_{{J_{\vec{\alpha}}}}\cap\Sigma_{{J_{\vec{\beta}}}}$, $\rho(A)=\emptyset$, $\forall{\vec{\beta}}\in\mathbb{S}^2$
(Theorem \ref{usa20}).
In the exceptional case when the Weyl function of $S$ is a constant, the set $\Xi_{{\vec{\alpha}}}$ increases considerably (Corollary \ref{usa70}).

The one-dimensional Schr\"{o}dinger differential expression with non-integrable singularity at zero (the limit-circle case at $x=0$)
is considered as an example of application (Proposition \ref{usa90}).

Throughout the paper, $\mathcal{D}(A)$ denotes the domain of a linear operator $A$. $A\upharpoonright_{\mathcal{D}}$ means the restriction of $A$ onto a set $\mathcal{D}$. The notation $\sigma(A)$ and $\rho(A)$ are used for the spectrum and the resolvent set of $A$.
The sign $\rule{2mm}{2mm}$ denotes the end of a proof.

\section{Sets $\Sigma_{{J_{\vec{\alpha}}}}$ and their properties}
\subsection{Preliminaries.}
Let $\mathfrak{H}$ be a Hilbert space with inner product
$(\cdot,\cdot)$ and let $J$ and $R$ be fundamental symmetries in $\mathfrak{H}$ satisfying (\ref{es545}).

Denote by ${\mathcal C}l_2(J,R):=\mbox{span}\{I, J, R, iJR\}$
a complex Clifford algebra with generating elements $J$ and $R$.
Since the operators $I, J, R$, and $iJR$ are linearly independent (due to (\ref{es545})), an arbitrary operator
$K\in{\mathcal C}l_2(J,R)$ can be presented as:
\begin{equation}\label{pp1}
K=\alpha_0{I}+\alpha_{1}J+\alpha_{2}R+\alpha_{3}iJR, \qquad \alpha_j\in\mathbb{C}
\end{equation}

\begin{lemma}[\cite{KT}]\label{ll1}
An operator $K$ defined by (\ref{pp1}) is a non-trivial fundamental symmetry in
$\mathfrak{H}$ (i.e., $K^2=I$, $K=K^*$, and
$K\not=I$) if and only if
\begin{equation}\label{es100}
K=\alpha_{1}J+\alpha_{2}R+\alpha_{3}iJR,
\end{equation}
where $\alpha_1^2+\alpha_2^2+\alpha_3^2=1$ and  $\alpha_j\in\mathbb{R}$.
\end{lemma}
\emph{Proof.} The reality of $\alpha_j$ in (\ref{es100}) follows from the self-adjointness of
$I, J, R,$ and $iJR$. The condition $\alpha_1^2+\alpha_2^2+\alpha_3^2=1$ is equivalent to the relation
$K^2=I$.
\rule{2mm}{2mm}

\begin{remark}\label{es42}
The formula (\ref{es100}) establishes a one-to-one correspondence between the set of non-trivial fundamental symmetries $K$ in
${\mathcal C}l_2(J,R)$ and vectors $\vec{\alpha}=(\alpha_1,\alpha_2,\alpha_3)$ of the unit sphere $\mathbb{S}^2$ in $\mathbb{R}^3$.
To underline this relationship we will use the notation $J_{\vec{\alpha}}$ for the fundamental symmetry $K$ determined
by (\ref{es100}) i.e.,
\begin{equation}\label{es100b}
J_{\vec{\alpha}}=\alpha_{1}J+\alpha_{2}R+\alpha_{3}iJR.
\end{equation}
In particular, this means that $J_{\vec{\alpha}}=J$ with $\vec{\alpha}=(1,0,0)$ and
$J_{\vec{\alpha}}=R$ when $\vec{\alpha}=(0,1,0)$.
\end{remark}

\begin{lemma}\label{ll2}
Let $\vec{\alpha}, \vec{\beta}\in{\mathbb{S}^2}$. Then:
 \begin{equation}\label{es6}
 J_{\vec{\alpha}}J_{\vec{\beta}}=-J_{\vec{\beta}}J_{\vec{\alpha}} \qquad \mbox{if and only if} \qquad \vec{\alpha}\cdot\vec{\beta}=0.
\end{equation}
and
\begin{equation}\label{es7}
 J_{\vec{\alpha}}+J_{\vec{\beta}}=|\vec{\alpha}+\vec{\beta}|J_{\frac{\vec{\alpha}+\vec{\beta}}{|\vec{\alpha}+\vec{\beta}|}} \qquad
 \mbox{if} \qquad \vec{\alpha}\not=-\vec{\beta}.
 \end{equation}
\end{lemma}
\emph{Proof.} It immediately follows from Lemma \ref{ll1} and identities (\ref{es545}), (\ref{es100b}).
\rule{2mm}{2mm}

\subsection{Definition and properties of $\Sigma_{{J_{\vec{\alpha}}}}$.}

\mbox{} \\

{\bf 1.}
Let $S$ be a closed symmetric operator \emph{with equal deficiency indices} in the Hilbert space $\mathfrak{H}$.
In what follows we suppose that $S$
\emph{commutes with all elements of} ${\mathcal C}l_2(J,R)$ or, that is equivalent, $S$ \emph{commutes with} $J$ and $R$:
\begin{equation}\label{es2}
SJ=JS, \qquad SR=RS
\end{equation}

Denote by $\Upsilon$ \emph{the set of all self-adjoint extensions $A$ of $S$ which commute with $J$ and $R$}:
\begin{equation}\label{es9}
\Upsilon=\{\ A\supset{S} \ : \ A^*=A, \quad AJ=JA, \quad AR=RA \ \}.
\end{equation}

It follows from (\ref{es100b}) and (\ref{es9}) that $\Upsilon$ contains self-adjoint extensions of $S$ \emph{which commute
with all fundamental symmetries}  $J_{\vec{\alpha}}\in{\mathcal C}l_2(J,R)$.

Let us fix one of them  $J_{\vec{\alpha}}$ and denote
by $(\mathfrak{H}, [\cdot,\cdot]_{J_{\vec{\alpha}}})$ the corresponding Krein space\footnote{We refer to \cite{AZ, GKS} for the terminology
of the Krein spaces theory} with
the indefinite inner product $[\cdot,\cdot]_{J_{\vec{\alpha}}}:=(J_{\vec{\alpha}}{\cdot}, \cdot)$.

Denote by $\Sigma_{J_{\vec{\alpha}}}$ \emph{the collection of all $J_{\vec{\alpha}}$-self-adjoint
extensions of $S$}:
\begin{equation}\label{es8}
\Sigma_{J_{\vec{\alpha}}}=\{\ A\supset{S} \ : \ J_{\vec{\alpha}}A^*=AJ_{\vec{\alpha}} \ \}.
\end{equation}

An operator $A\in\Sigma_{J_{\vec{\alpha}}}$ is a self-adjoint extension of $S$ with respect to the indefinite metric
$[\cdot,\cdot]_{J_{\vec{\alpha}}}$.

\begin{proposition}\label{um1}
The following relation holds
$$
\bigcap_{\forall\vec{\alpha}\in{\mathbb{S}^2}}\Sigma_{J_{\vec{\alpha}}}=\Upsilon.
$$
\end{proposition}
\emph{Proof.}
It follows from the definitions above that $\Sigma_{J_{\vec{\alpha}}}\supset\Upsilon$. Therefore,
 $$
 \bigcap_{\forall\vec{\alpha}\in{\mathbb{S}^2}}\Sigma_{J_{\vec{\alpha}}}\supset\Upsilon.
 $$

Let $A\in\bigcap_{\forall\vec{\alpha}\in{\mathbb{S}^2}}\Sigma_{J_{\vec{\alpha}}}$.
In particular, this means that $A\in\Sigma_{J}$, $A\in\Sigma_{R}$, and $A\in\Sigma_{iJR}$. It follows from the first two relations that
$JA^*=AJ$ and $RA^*=AR$. Therefore, $iJRA^*=iJAR=A^*iJR$.
Simultaneously, $iJRA^*=AiJR$ since $A\in\Sigma_{iJR}$.
Comparing the obtained relations we deduce that $A^*iJR=AiJR$  and hence, $A^*=A$.
Thus $A$ is a self-adjoint operator and it commutes with
an arbitrary fundamental symmetry $J_{\vec{\alpha}}\in{\mathcal C}l_2(J,R)$. Therefore,
$A\in\Upsilon$. Proposition \ref{um1} is proved.
\rule{2mm}{2mm}

\smallskip

Simple analysis of the proof of Proposition \ref{um1} leads to the conclusion that
$$
\Sigma_{J_{\vec{\alpha}}}\cap\Sigma_{J_{\vec{\beta}}}\cap\Sigma_{J_{\vec{\gamma}}}=\Upsilon
$$
\emph{for any three linearly independent vectors} $\vec{\alpha}, \vec{\beta}, \vec{\gamma}\in{\mathbb{S}^2}$.
However,
\begin{equation}\label{es444}
\Sigma_{J_{\vec{\alpha}}}\cap\Sigma_{J_{\vec{\beta}}}\supset\Upsilon, \qquad \forall{\vec{\alpha}, \vec{\beta}\in{\mathbb{S}^2}}
\end{equation}
and the intersection $\Sigma_{J_{\vec{\alpha}}}\cap\Sigma_{J_{\vec{\beta}}}$ contains operators $A$ \emph{with empty resolvent set}
 (i.e., $\rho(A)=\emptyset$ or, that is equivalent, $\sigma(A)=\mathbb{C}$).
Let us discuss this phenomena in detail.

Consider two linearly independent vectors $\vec{\alpha}, \vec{\beta}\in{\mathbb{S}^2}$.
If $\vec{\alpha}\cdot\vec{\beta}\not=0$, we define
new vector $\vec{\beta'}$ in $\mathbb{S}^2$:
$$
\vec{\beta'}=\frac{\vec{\alpha}+c\vec{\beta}}{|\vec{\alpha}+c\vec{\beta}|}, \qquad c=-\frac{1}{\vec{\alpha}\cdot\vec{\beta}}
$$
such that $\vec{\alpha}\cdot\vec{\beta'}=0$. Then the fundamental symmetry
\begin{equation}\label{es89}
J_{\vec{\beta'}}=\frac{1}{|\vec{\alpha}+c\vec{\beta}|}J_{\vec{\alpha}}+\frac{c}{|\vec{\alpha}+c\vec{\beta}|}J_{\vec{\beta}}
\end{equation}
anti-commutes with $J_{\vec{\alpha}}$ (due to Lemma \ref{ll2}).

The operator
\begin{equation}\label{new78}
J_{\vec{\gamma}}=iJ_{\vec{\alpha}}J_{\vec{\beta'}}=\left|\begin{array}{ccc}
J & R & iJR \\
\alpha_1 & \alpha_2 & \alpha_3 \\
\beta'_1 & \beta'_2 & \beta'_3
\end{array}\right|
\end{equation}
is a fundamental symmetry in $\mathfrak{H}$ which commutes with $S$. Therefore,
the orthogonal decomposition of $\mathfrak{H}$ constructed by $J_{\vec{\gamma}}$:
\begin{equation}\label{es90}
\mathfrak{H}=\mathfrak{H}_{+}^{\gamma}\oplus\mathfrak{H}^{\gamma}_{-}, \qquad  \mathfrak{H}_{+}^{\gamma}=\frac{1}{2}(I+J_{\vec{\gamma}})\mathfrak{H}, \quad
\mathfrak{H}_{-}^{\gamma}=\frac{1}{2}(I-J_{\vec{\gamma}})\mathfrak{H}
\end{equation}
reduces $S$:
\begin{equation}\label{es1b}
S=\left(\begin{array}{cc}
S_{\gamma+} & 0 \\
0 & S_{\gamma-} \end{array}\right), \qquad S_{\gamma+}=S\upharpoonright_{\mathfrak{H}_{+}^{\gamma}}, \quad S_{\gamma-}=S\upharpoonright_{\mathfrak{H}_{-}^{\gamma}}.
\end{equation}

Since $J_{\vec{\gamma}}$ anti-commutes with
$J_{\vec{\alpha}}$ (see (\ref{new78})), the operator $J_{\vec{\alpha}}$ maps $\mathfrak{H}_{\pm}^{\gamma}$
onto $\mathfrak{H}_{\mp}^{\gamma}$ and operators $S_{\gamma+}$ and $S_{\gamma-}$ are unitarily equivalent.
Precisely, $S_{\gamma-}=J_{\vec{\alpha}}S_{\gamma+}J_{\vec{\alpha}}$.
This means that $S_{\gamma+}$ and $S_{\gamma-}$ have equal deficiency indices.\footnote{this also implies that
the symmetric operator $S$ commuting with ${\mathcal C}l_2(J,R)$ may have only \emph{even} deficiency indices.}

Denote
\begin{equation}\label{ew1}
A_\gamma=\left(\begin{array}{cc}
S_{\gamma+} & 0 \\
0 & S_{\gamma-}^* \end{array}\right), \qquad  A_\gamma^*=\left(\begin{array}{cc}
S_{\gamma+}^* & 0 \\
0 & S_{\gamma-} \end{array}\right)
\end{equation}
The operators $A_\gamma$ and $A_\gamma^*$ are extensions of $S$ and $\sigma(A_\gamma)=\sigma(A_\gamma^*)=\mathbb{C}$
(since $S_{\gamma\pm}$ are symmetric operators), i.e.,
these operators \emph{have empty resolvent set}.

\begin{theorem}\label{um2}
Let  $\vec{\alpha}, \vec{\beta}\in{\mathbb{S}^2}$ be
linearly independent vectors. Then the operators $A_\gamma$ and $A_\gamma^*$ defined by (\ref{ew1}) belong to
$\Sigma_{J_{\vec{\alpha}}}\cap\Sigma_{J_{\vec{\beta}}}$.
\end{theorem}
\emph{Proof.} Assume that $A\in\Sigma_{J_{\vec{\alpha}}}\cap\Sigma_{J_{\vec{\beta}}}$, where $\vec{\alpha}, \vec{\beta}\in{\mathbb{S}^2}$ are linearly
independent vectors. Then
\begin{equation}\label{esse34}
J_{\vec{\alpha}}A=A^*J_{\vec{\alpha}}, \qquad  J_{\vec{\beta}}A=A^*J_{\vec{\beta}}
\end{equation}
and hence, $J_{\vec{\beta'}}A=A^*J_{\vec{\beta'}}$ due to (\ref{es89}). In that case
$$
J_{\vec{\gamma}}A=iJ_{\vec{\alpha}}J_{\vec{\beta'}}A=iJ_{\vec{\alpha}}A^*J_{\vec{\beta'}}=iAJ_{\vec{\alpha}}J_{\vec{\beta'}}=AJ_{\vec{\gamma}},
$$
where $J_{\vec{\gamma}}=iJ_{\vec{\alpha}}J_{\vec{\beta'}}$ (see (\ref{new78})) is the fundamental symmetry in $\mathfrak{H}$.

Since $A$ commutes with $J_{\vec{\gamma}}$, the decomposition (\ref{es90}) reduces $A$ and
\begin{equation}\label{esse5}
A=\left(\begin{array}{cc}
A_{+} & 0 \\
0 & A_{-} \end{array}\right), \qquad A_{+}=A\upharpoonright_{\mathfrak{H}_{+}^{\gamma}}, \quad A_{-}=A\upharpoonright_{\mathfrak{H}_{-}^{\gamma}},
\end{equation}
where $S_{\gamma+}{\subseteq}A_{+}{\subseteq}S_{\gamma+}^*$ and $S_{\gamma-}{\subseteq}A_{-}{\subseteq}S_{\gamma-}^*$.
This means that
\begin{equation}\label{usa1}
A^*=\left(\begin{array}{cc}
A_{+}^* & 0 \\
0 & A_{-}^* \end{array}\right), \qquad A_{+}^*=A^*\upharpoonright_{\mathfrak{H}_{+}^{\gamma}}, \quad A_{-}^*=A^*\upharpoonright_{\mathfrak{H}_{-}^{\gamma}},
\end{equation}

Since $J_{\vec{\alpha}}$ anti-commutes with $J_{\vec{\gamma}}$, the operator $J_{\vec{\alpha}}$ maps $\mathfrak{H}_{\pm}^{\gamma}$ onto $\mathfrak{H}_{\mp}^{\gamma}$. Therefore, the first relation in (\ref{esse34})  can be
rewritten with the use of formulas (\ref{esse5}) and (\ref{usa1}) as follows:
\begin{equation}\label{esse6}
J_{\vec{\alpha}}Ax=J_{\vec{\alpha}}(A_{+}x_++A_{-}x_-)=A_{-}^*J_{\vec{\alpha}}x_++A_{+}^*J_{\vec{\alpha}}x_-=A^*J_{\vec{\alpha}}x,
\end{equation}
where $x=x_++x_-\in\mathcal{D}(A), \ x_\pm\in\mathcal{D}(A_\pm)$.

The identity (\ref{esse6}) holds for all $x_\pm\in\mathcal{D}(A_\pm)$. This means that
\begin{equation}\label{esse7}
J_{\vec{\alpha}}A_{+}=A_{-}^*J_{\vec{\alpha}}, \qquad J_{\vec{\alpha}}A_{-}=A_{+}^*J_{\vec{\alpha}}
\end{equation}

It follows from (\ref{es89}) and (\ref{new78}) that the fundamental symmetry $J_{\vec{\beta}}$ anti-commutes with $J_{\vec{\gamma}}$.
Repeating the arguments above for the second relation in (\ref{esse34}) we obtain
\begin{equation}\label{esse8}
J_{\vec{\beta}}A_{+}=A_{-}^*J_{\vec{\beta}}, \qquad J_{\vec{\beta}}A_{-}=A_{+}^*J_{\vec{\beta}}.
\end{equation}

Thus an operator $A$ belongs to $\Sigma_{J_{\vec{\alpha}}}\cap\Sigma_{J_{\vec{\beta}}}$
\emph{if and only if its counterparts $A_+$ and $A_-$ in (\ref{esse5}) satisfy relations (\ref{esse7}) and (\ref{esse8})}.
In particular, these relations are satisfied for the cases when $A_+=S_{\gamma+}, \ A_-=S_{\gamma-}^*$ and
$A_+=S_{\gamma+}^*, \ A_-=S_{\gamma-}$. Hence, the
operators $A_\gamma$, $A_\gamma^*$  defined by (\ref{ew1})  belong to $\Sigma_{J_{\vec{\alpha}}}\cap\Sigma_{J_{\vec{\beta}}}$.
Theorem \ref{um2} is proved.
\rule{2mm}{2mm}
\begin{remark}\label{esse60} The operators $A_\gamma$ and $A_\gamma^*$ constructed above
depend on the choice of $\vec{\beta}\in\mathbb{S}^2$.
Considering various vectors $\vec{\beta}\in\mathbb{S}^2$ in (\ref{es89}), (\ref{new78}), we obtain a collection of fundamental symmetries $J_{\vec{\gamma}(\vec{\beta})}$. This gives rise to a one-parameter set of different operators $A_{\gamma(\vec{\beta})}$ and $A_{\gamma(\vec{\beta})}^*$ with empty resolvent set which belong to $\Sigma_{J_{\vec{\alpha}}}$.
\end{remark}

\begin{corollary}\label{usa3}
Let  $\vec{\alpha}, \vec{\beta}\in{\mathbb{S}^2}$ be
linearly independent vectors and let (\ref{es90}) be the decomposition of $\mathfrak{H}$ constructed by these vectors.
Then, with respect to (\ref{es90}), all operators $A\in\Sigma_{J_{\vec{\alpha}}}\cap\Sigma_{J_{\vec{\beta}}}$ are described by the formula
\begin{equation}\label{usa8}
A=\left(\begin{array}{cc}
A_{+} & 0 \\
0 & J_{\vec{\alpha}}A_{+}^*J_{\vec{\alpha}} \end{array}\right),
\end{equation}
where $A_+$ is an arbitrary intermediate extension of $S_{\gamma+}=S\upharpoonright_{\mathfrak{H}_{+}^{\gamma}}$ (i.e., $S_{\gamma+}{\subseteq}A_{+}{\subseteq}S_{\gamma+}^*$)
\end{corollary}
\emph{Proof.} If $A\in\Sigma_{J_{\vec{\alpha}}}\cap\Sigma_{J_{\vec{\beta}}}$, then the presentation (\ref{usa8}) follows from
(\ref{esse5}) and the second identity in (\ref{esse7}).

Conversely, assume that an operator $A$ is defined by (\ref{usa8}). Since $J_{\vec{\alpha}}$
and $J_{\vec{\beta}}$  anti-commute with $J_{\vec{\gamma}}$, they admit the presentations
$J_{\vec{\alpha}}=\left(\begin{array}{cc}
0 & J_{\vec{\alpha}} \\
J_{\vec{\alpha}} & 0 \end{array}\right)$ and $J_{\vec{\beta}}=\left(\begin{array}{cc}
0 & J_{\vec{\beta}} \\
J_{\vec{\beta}} & 0 \end{array}\right)$ with respect to (\ref{es90}). Then, the operator equality
$J_{\vec{\alpha}}A=A^*J_{\vec{\alpha}}$ is established by the direct multiplication of the corresponding operator entries.
The same procedure for $J_{\vec{\beta}}A=A^*J_{\vec{\beta}}$ leads to the verification of relations
\begin{equation}\label{usa19}
J_{\vec{\alpha}}J_{\vec{\beta}}A_+=A_+J_{\vec{\alpha}}J_{\vec{\beta}}, \qquad J_{\vec{\beta}}J_{\vec{\alpha}}A_+^*=A_+^*J_{\vec{\beta}}J_{\vec{\alpha}}.
\end{equation}
To this end we recall that $J_{\vec{\gamma}}$ commutes with $A_+$ and
$$
J_{\vec{\alpha}}J_{\vec{\beta}}=-\frac{1}{c}I-i\frac{|\vec{\alpha}+c\vec{\beta}|}{c}J_{\vec{\gamma}}
$$
due to (\ref{es89}) and (\ref{new78}). Therefore, $A_+$ commutes with $J_{\vec{\alpha}}J_{\vec{\beta}}$ and
the first relation in (\ref{usa19}) holds. The second relation is established in the same manner, if we take into
account that $J_{\vec{\gamma}}$ commutes with $A_+^*$ and $J_{\vec{\beta}}J_{\vec{\alpha}}=-\frac{1}{c}I+i\frac{|\vec{\alpha}+c\vec{\beta}|}{c}J_{\vec{\gamma}}$.
Corollary \ref{usa3} is proved. \rule{2mm}{2mm}
\begin{remark}
It follows from the proof that the choice of $J_{\vec{\alpha}}$ in (\ref{usa8}) is not essential and
the similar description of $\Sigma_{J_{\vec{\alpha}}}\cap\Sigma_{J_{\vec{\beta}}}$ can be obtained with the help
of $J_{\vec{\beta}}$.
\end{remark}

\smallskip

{\bf 2.} Denote
\begin{equation}\label{es10}
W_{\vec{\alpha},\vec{\beta}}=\left\{\begin{array}{l}
\displaystyle{J_{\frac{\vec{\alpha}+\vec{\beta}}{|\vec{\alpha}+\vec{\beta}|}}} \quad \mbox{if} \quad \vec{\alpha}\not=-\vec{\beta}; \vspace{3mm} \\
I \qquad \mbox{if} \qquad \vec{\alpha}=-\vec{\beta}
\end{array}\right.
\end{equation}

It is clear that $W_{\vec{\alpha},\vec{\beta}}$ is a fundamental symmetry in $\mathfrak{H}$
and $W_{\vec{\alpha},\vec{\beta}}=W_{\vec{\beta}, \vec{\alpha}}$ for any $\vec{\alpha},\vec{\beta}\in{\mathbb{S}^2}$.

\begin{theorem}\label{es12a}
For any $\vec{\alpha},\vec{\beta}\in{\mathbb{S}^2}$ the sets $\Sigma_{J_{\vec{\alpha}}}$ and $\Sigma_{J_{\vec{\beta}}}$ are unitarily equivalent and $A\in\Sigma_{J_{\vec{\alpha}}}$ if and only if $W_{\vec{\alpha},\vec{\beta}}AW_{\vec{\alpha},\vec{\beta}}\in\Sigma_{J_{\vec{\beta}}}$.
\end{theorem}
\emph{Proof.} Since
$J_{-\vec{\alpha}}=-J_{\vec{\alpha}}$ (see (\ref{es100b})), the sets $\Sigma_{J_{\vec{\alpha}}}$ and $\Sigma_{J_{-\vec{\alpha}}}$ coincide and therefore, the case $\vec{\alpha}=-\vec{\beta}$ is trivial.

Assume that $A\in\Sigma_{J_{\vec{\alpha}}}$, $\vec{\alpha}\not=-\vec{\beta}$   and consider the operator
$W_{\vec{\alpha},\vec{\beta}}AW_{\vec{\alpha},\vec{\beta}}$, which we denote $B=W_{\vec{\alpha},\vec{\beta}}AW_{\vec{\alpha},\vec{\beta}}$ for brevity.
Taking into account that $S$ commutes with
$J_{\vec{\alpha}}$ for any choice of $\vec{\alpha}\in{\mathbb{S}^2}$, we deduce from (\ref{es10}) that $
W_{\vec{\alpha},\vec{\beta}}S=SW_{\vec{\alpha},\vec{\beta}} \quad \mbox{and} \quad W_{\vec{\alpha},\vec{\beta}}S^*=S^*W_{\vec{\alpha},\vec{\beta}}.$ This means that
$Bx=W_{\vec{\alpha},\vec{\beta}}AW_{\vec{\alpha},\vec{\beta}}x=W_{\vec{\alpha},\vec{\beta}}SW_{\vec{\alpha},\vec{\beta}}x=Sx$ for all
$x\in\mathcal{D}(S)$ and $By=W_{\vec{\alpha},\vec{\beta}}AW_{\vec{\alpha},\vec{\beta}}y=W_{\vec{\alpha},\vec{\beta}}S^*W_{\vec{\alpha},\vec{\beta}}y=S^*y$ for all
$y\in\mathcal{D}(B)=W_{\vec{\alpha},\vec{\beta}}\mathcal{D}(A)$. Therefore, $B$ is an intermediate extension of $S$ (i.e. $S\subset{B}\subset{S^*}$).

It follows from (\ref{es7}) and (\ref{es10}) that
\begin{equation}\label{es15}
J_{\vec{\beta}}W_{\vec{\alpha},\vec{\beta}}=J_{\vec{\beta}}\frac{J_{\vec{\alpha}}+J_{\vec{\beta}}}{|\vec{\alpha}+\vec{\beta}|}=\frac{J_{\vec{\beta}}J_{\vec{\alpha}}+I}{|\vec{\alpha}+\vec{\beta}|}=\frac{J_{\vec{\beta}}+J_{\vec{\alpha}}}{|\vec{\alpha}+\vec{\beta}|}J_{\vec{\alpha}}=W_{\vec{\alpha},\vec{\beta}}J_{\vec{\alpha}}.
\end{equation}
Using (\ref{esse34}) and (\ref{es15}), we arrive at the conclusion that
$$
J_{\vec{\beta}}B^*=J_{\vec{\beta}}W_{\vec{\alpha},\vec{\beta}}A^*W_{\vec{\alpha},\vec{\beta}}=W_{\vec{\alpha},\vec{\beta}}AJ_{\vec{\alpha}}W_{\vec{\alpha},\vec{\beta}}=W_{\vec{\alpha},\vec{\beta}}AW_{\vec{\alpha},\vec{\beta}}J_{\vec{\beta}}=BJ_{\vec{\beta}}.
$$
Therefore, condition $A\in\Sigma_{J_{\vec{\alpha}}}$ implies that ${B=W_{\vec{\alpha},\vec{\beta}}AW_{\vec{\alpha},\vec{\beta}}}\in\Sigma_{J_{\vec{\beta}}}$.
The inverse implication  $B\in\Sigma_{J_{\vec{\beta}}}\Rightarrow{A=W_{\vec{\alpha},\vec{\beta}}BW_{\vec{\alpha},\vec{\beta}}}\in\Sigma_{J_{\vec{\alpha}}}$ is established in the same manner. Theorem \ref{es12a} is proved.
\rule{2mm}{2mm}
\begin{remark}
Due to Theorem \ref{es12a}, for any $J_{\vec{\alpha}}$-self-adjoint extension $A\in\Sigma_{J_{\vec{\alpha}}}$
there exists a unitarily equivalent $J_{\vec{\beta}}$-self-adjoint extension $B\in\Sigma_{J_{\vec{\beta}}}$.  This means that, the spectral analysis of operators from $\bigcup_{\alpha\in{\mathbb{S}^2}}\Sigma_{J_{\vec{\alpha}}}$ can be reduced
to the spectral analysis of $J_{\vec{\alpha}}$-self-adjoint extensions from  $\Sigma_{J_{\vec{\alpha}}}$,
where $\vec{\alpha}$ is a fixed vector from $\mathbb{S}^2$.
\end{remark}

\subsection{Boundary triplets and Weyl function.}

\mbox{} \\

{\bf 1.} Let $S$ be a closed symmetric operator with equal deficiency indices in the Hilbert space $\mathfrak{H}$.
A triplet $(\mathcal{H}, \Gamma_0, \Gamma_1)$, where $\mathcal{H}$ is an
auxiliary Hilbert space and $\Gamma_0$, $\Gamma_1$ are linear
mappings of $\mathcal{D}(S^*)$ into $\mathcal{H}$, is called a
\emph{boundary triplet of} $S^*$ if the abstract Green identity
\begin{equation}\label{new2}
(S^*x, y)-(x, S^*y)=(\Gamma_1x, \Gamma_0y)_{\mathcal{H}}-(\Gamma_0x,
\Gamma_1y)_{\mathcal{H}}, \quad  x, y\in\mathcal{D}(S^*)
\end{equation}
is satisfied and the map $(\Gamma_0,
\Gamma_1):\mathcal{D}(S^*)\to\mathcal{H}\oplus\mathcal{H}$ is
surjective \cite{Gor}.
\begin{lemma}\label{es16}
Assume that $S$ satisfies the commutation relations (\ref{es2}) and
$J_{\vec{\tau}}, J_{\vec{\gamma}}\in{\mathcal Cl}_2(J,R)$ are fixed anti-commuting fundamental symmetries.
Then there exists a boundary triplet
 $(\mathcal{H}, \Gamma_0, \Gamma_1)$ of $S^*$ such that
the formulas
\begin{equation}\label{ea6a}
{\mathcal J}_{\vec{\tau}}\Gamma_j:=\Gamma_j{J}_{\vec{\tau}}, \qquad
{\mathcal J}_{\vec{\gamma}}\Gamma_j:=\Gamma_j{J}_{\vec{\gamma}}, \qquad j=0,1
\end{equation}
correctly define anti-commuting fundamental symmetries ${\mathcal J}_{\vec{\tau}}$ and ${\mathcal J}_{\vec{\gamma}}$ in the Hilbert space $\mathcal{H}$.
\end{lemma}
\emph{Proof.}
If $S$ satisfies (\ref{es2}), then $S$
commutes with an arbitrary fundamental symmetry $J_{\vec{\gamma}}\in{\mathcal Cl}_2(J,R)$ and hence,
$S$ admits the representation (\ref{es1b}) for any vector $\vec{\gamma}\in\mathbb{S}^2$.

Let $S_{\gamma+}$ be a symmetric operator in $\mathfrak{H}_+^\gamma$ from (\ref{es1b})
and let $(N, \Gamma_0^+, \Gamma_1^+)$ be an arbitrary
 boundary triplet of $S^*_{\gamma+}$.

Since $J_{\vec{\tau}}$ anti-commutes with $J_{\vec{\gamma}}$, the symmetric operator $S_{\gamma-}$ in (\ref{es1b})
can be described as $S_{\gamma-}=J_{\vec{\tau}}S_{\gamma+}J_{\vec{\tau}}$.
This means that $(N, \Gamma_0^+{J_{\vec{\tau}}},
\Gamma_1^+{J_{\vec{\tau}}})$ is a boundary triplet of $S_{\gamma-}^*$.

It is easy to see that the operators
\begin{equation}\label{esse41}
\Gamma_jf=\Gamma_j(f_++f_-)=\left(\begin{array}{c}
\Gamma_j^+f_+ \\
\Gamma_j^+J_{\vec{\tau}}f_-
\end{array}\right)
\end{equation}
($f=f_++f_-\in\mathcal{D}(S^*), \ f_\pm\in\mathcal{D}(S_{\gamma\pm}^*)$) map $\mathcal{D}(S^*)$ onto
the Hilbert space
$$
\mathcal{H}=\mathcal{H}_+\oplus\mathcal{H}_-, \qquad  \mathcal{H}_+=\left(\begin{array}{c}
N \\
0
\end{array}\right), \quad \mathcal{H}_-=\left(\begin{array}{c}
0 \\
N
\end{array}\right)
$$
and they form a boundary triplet $(\mathcal{H}, \Gamma_0, \Gamma_1)$ of $S^*$ which
satisfies (\ref{ea6a}) with
\begin{equation}\label{usa31}
\mathcal{J}_{\vec{\tau}}=\left(\begin{array}{cc}
0 &   I\\
I & 0
\end{array}\right), \qquad \mathcal{J}_{\vec{\gamma}}=\left(\begin{array}{cc}
I & 0  \\
0 & -I
\end{array}\right).
\end{equation}
It is clear that $\mathcal{J}_{\vec{\tau}}$ and $\mathcal{J}_{\vec{\gamma}}$ are anti-commuting
fundamental symmetries in the Hilbert space $\mathcal{H}$.
\rule{2mm}{2mm}

\begin{remark}\label{es21}
Since $J$ and $R$ can be expressed as linear combinations of ${J}_{\vec{\tau}}$, ${J}_{\vec{\gamma}}$, and
$i{J}_{\vec{\tau}}{J}_{\vec{\gamma}}$, formulas (\ref{ea6a}) imply that
$$
{\mathcal J}\Gamma_j:=\Gamma_j{J}, \qquad
{\mathcal R}\Gamma_j:=\Gamma_j{R}, \qquad j=0,1,
$$
where ${\mathcal J}$ and ${\mathcal R}$ are
anti-commuting fundamental symmetries in $\mathcal{H}$.
Therefore, an arbitrary boundary triplet
 $(\mathcal{H}, \Gamma_0, \Gamma_1)$ of $S^*$ with property (\ref{ea6a})
 allows one \emph{to establish a bijective correspondence between
elements of the initial Clifford algebra ${\mathcal C}l_2({J}, R)$ and its
'image' ${\mathcal C}l_2(\mathcal{J}, \mathcal{R})$ in the auxiliary space $\mathcal{H}$}.
In particular, \emph{for every} ${J}_{\vec{\alpha}}\in{\mathcal C}l_2({J},{R})$ defined by (\ref{es100b}),
\begin{equation}\label{es45}
{\mathcal {J}_{\vec{\alpha}}}\Gamma_j=\Gamma_j{{J}_{\vec{\alpha}}},  \qquad  j=0,1,
\end{equation}
where
$
\mathcal{{J}_{\vec{\alpha}}}=\alpha_{1}\mathcal{J}+\alpha_{2}\mathcal{R}+\alpha_{3}i\mathcal{JR}
$
belongs to ${\mathcal C}l_2({\mathcal J},{\mathcal R})$.
\end{remark}

\smallskip

{\bf 2.} Let $(\mathcal{H}, \Gamma_0,
\Gamma_1)$ be a boundary triplet of $S^*$.
The Weyl function of $S$ associated with $(\mathcal{H}, \Gamma_0,
\Gamma_1)$ is defined as follows:
\begin{equation}\label{neww65}
M(\mu)\Gamma_0f_\mu=\Gamma_1f_\mu, \quad
\forall{f}_\mu\in\ker(S^*-\mu{I}), \quad
\forall\mu\in\mathbb{C}\setminus\mathbb{R}.
\end{equation}

\begin{lemma}\label{es30}
Let $(\mathcal{H}, \Gamma_0, \Gamma_1)$ be a boundary triplet
of $S^*$ with properties (\ref{ea6a}). Then the corresponding Weyl function $M(\cdot)$ commutes with every
fundamental symmetry ${\mathcal J}_{\vec{\alpha}}\in{\mathcal C}l_2({\mathcal J},{\mathcal R})$:
$$
M(\mu){\mathcal J}_{\vec{\alpha}}={\mathcal J}_{\vec{\alpha}}M(\mu), \qquad \forall\mu\in\mathbb{C}\setminus\mathbb{R}.
$$
\end{lemma}
\emph{Proof.} It follows from (\ref{es100b}) and (\ref{es2}) that
$S^*J_{\vec{\alpha}}=J_{\vec{\alpha}}S^*$ for all $\vec{\alpha}\in\mathbb{S}^2$. Therefore,
${J}_{\vec{\alpha}}:\ker(S^*-\mu{I})\to\ker(S^*-\mu{I})$. In that case, relations (\ref{es45}) and (\ref{neww65}) lead to
$M(\mu){\mathcal J}_{\vec{\alpha}}\Gamma_0f_\mu={\mathcal J}_{\vec{\alpha}}\Gamma_1f_\mu$. Thus,
${\mathcal J}_{\vec{\alpha}}M(\mu){\mathcal J}_{\vec{\alpha}}=M(\mu)$ or $M(\mu){\mathcal J}_{\vec{\alpha}}={\mathcal J}_{\vec{\alpha}}M(\mu)$. \rule{2mm}{2mm}

\subsection{Description of $\Sigma_{J_{\vec{\alpha}}}$ in terms of boundary triplets.}
\begin{theorem}\label{new45}
Let $(\mathcal{H}, \Gamma_0, \Gamma_1)$ be a boundary triplet
of $S^*$ with properties (\ref{ea6a}) for a fixed anti-commuting fundamental symmetries $J_{\vec{\tau}}, J_{\vec{\gamma}}\in{\mathcal Cl}_2(J,R)$ and let $J_{\vec{\alpha}}$ be an arbitrary fundamental symmetry from ${\mathcal Cl}_2(J,R)$. Then operators $A\in\Sigma_{J_{\vec{\alpha}}}$ coincide with the restriction of $S^*$ onto the domains
\begin{equation}\label{sas98}
\mathcal{D}(A)=\{f\in\mathcal{D}(S^*) \ : \ U({{\mathcal J}_{\vec{\alpha}}}\Gamma_1+i\Gamma_0)f=({{\mathcal J}_{\vec{\alpha}}\Gamma_1-i\Gamma_0)f \}},
\end{equation}
where $U$ runs the set of unitary operators in $\mathcal{H}$.

The correspondence $A\leftrightarrow{U}$ determined by (\ref{sas98}) is a bijection between the set
$\Sigma_{J_{\vec{\alpha}}}$ of all $J_{\vec{\alpha}}$-self-adjoint extensions of $S$ and the set of unitary operators in $\mathcal{H}$.
\end{theorem}

{\it Proof.} An operator $A$ is $J_{\vec{\alpha}}$-self-adjoint extension of $S$ if and only if
$J_{\vec{\alpha}}A$ is a self-adjoint extension of the symmetric operator $J_{\vec{\alpha}}S$.
Since
\begin{equation}\label{neww12}
(J_{\vec{\alpha}}S)^*=S^*J_{\vec{\alpha}}=J_{\vec{\alpha}}S^*,
\end{equation}
the Green identity (\ref{new2}) can be rewritten with the use of (\ref{es45})
as follows:
$$
(S^*J_{\vec{\alpha}}x, y)-(x, S^*J_{\vec{\alpha}}y)=({\mathcal J}_{\vec{\alpha}}\Gamma_1x, \Gamma_0y)_{\mathcal{H}}-(\Gamma_0x,
{\mathcal J}_{\vec{\alpha}}\Gamma_1y)_{\mathcal{H}}.
$$
Recalling the definition of boundary triplet we conclude that $(\mathcal{H}, \Gamma_0, {\mathcal J}_{\vec{\alpha}}\Gamma_1)$ is a boundary triplet of
$J_{\vec{\alpha}}S$.
Therefore \cite[Chapter 3, Theorem 1.6]{Gor}, self-adjoint extensions $J_{\vec{\alpha}}A$ of $J_{\vec{\alpha}}S$ coincide with the restriction
of $(J_{\vec{\alpha}}S)^*$ onto
$$
\mathcal{D}(J_{\vec{\alpha}}A)=\{f\in\mathcal{D}((J_{\vec{\alpha}}S)^*) \ : \ U({\mathcal J}_{\vec{\alpha}}\Gamma_1+i\Gamma_0)f=({\mathcal J}_{\vec{\alpha}}\Gamma_1-i\Gamma_0)f \}
$$
where $U$ runs the set of unitary operators in $\mathcal{H}$ and the correspondence $J_{\vec{\alpha}}A\leftrightarrow{U}$ is bijective.
 By virtue of (\ref{neww12}), $\mathcal{D}((J_{\vec{\alpha}}S)^*)=\mathcal{D}(S^*)$. Hence, $J_{\vec{\alpha}}$-self-adjoint extensions $A$ of $S$ coincide with the restriction of $S^*$ onto
$\mathcal{D}(J_{\vec{\alpha}}A)$ that implies (\ref{sas98}).  \rule{2mm}{2mm}

\begin{corollary}\label{esse50}
If $A\in\Sigma_{J_{\vec{\alpha}}}$ and $A\leftrightarrow{U}$ in (\ref{sas98}), then the
$J_{\vec{\beta}}$-self-adjoint operator ${B=W_{\vec{\alpha},\vec{\beta}}AW_{\vec{\alpha},\vec{\beta}}}\in\Sigma_{J_{\vec{\beta}}}$
($\vec{\alpha}\not=-\vec{\beta}$) is determined by the formula
$$
B=S^*\upharpoonright\{g\in\mathcal{D}(S^*) \ : \ {\mathcal W}_{\vec{\alpha},\vec{\beta}}U{\mathcal W}_{\vec{\alpha},\vec{\beta}}({\mathcal J}_{\vec{\beta}}\Gamma_1+i\Gamma_0)g=({\mathcal J}_{\vec{\beta}}\Gamma_1-i\Gamma_0)g \},
$$
where ${\mathcal W}_{\vec{\alpha},\vec{\beta}}={{\mathcal J}_{\frac{\vec{\alpha}+\vec{\beta}}{|\vec{\alpha}+\vec{\beta}|}}}$
is a fundamental symmetry in $\mathcal{H}$.
\end{corollary}
\emph{Proof.} Let  $A\in\Sigma_{J_{\vec{\alpha}}}$. Then ${B=W_{\vec{\alpha},\vec{\beta}}AW_{\vec{\alpha},\vec{\beta}}}\in\Sigma_{J_{\vec{\beta}}}$ by Theorem \ref{es12a} and, in view of Theorem \ref{new45},
\begin{equation}\label{sas100}
B=S^*\upharpoonright\{g\in\mathcal{D}(S^*) \ : \ U'({\mathcal J}_{\vec{\beta}}\Gamma_1+i\Gamma_0)g=({\mathcal J}_{\vec{\beta}}\Gamma_1-i\Gamma_0)g \},
\end{equation}
where $U'$ is a unitary operator in $\mathcal{H}$.

It follows from the definition of $B$ that $f\in\mathcal{D}(A)$ if and only if $g=W_{\vec{\alpha},\vec{\beta}}f\in\mathcal{D}(B)$.
Hence, we can rewrite (\ref{sas100}) with the use of (\ref{es15}):
\begin{equation}\label{esse51}
U'({\mathcal J}_{\vec{\beta}}\Gamma_1+i\Gamma_0)g=U'\mathcal{W}_{\vec{\alpha},\vec{\beta}}({\mathcal J}_{\vec{\alpha}}\Gamma_1+i\Gamma_0)f=({\mathcal J}_{\vec{\beta}}\Gamma_1-i\Gamma_0)g=\mathcal{W}_{\vec{\alpha},\vec{\beta}}({\mathcal J}_{\vec{\alpha}}\Gamma_1-i\Gamma_0)f,
\end{equation}
where $\mathcal{W}_{\vec{\alpha},\vec{\beta}}\Gamma_j=\Gamma_j{W}_{\vec{\alpha},\vec{\beta}}, \ j=0,1$ (cf. (\ref{es45})).

It follows from (\ref{es10}) that ${\mathcal W}_{\vec{\alpha},\vec{\beta}}={{\mathcal J}_{\frac{\vec{\alpha}+\vec{\beta}}{|\vec{\alpha}+\vec{\beta}|}}}$ and hence, $\mathcal{W}_{\vec{\alpha},\vec{\beta}}$ is
a fundamental symmetry in $\mathcal{H}$. Comparing (\ref{esse51}) with (\ref{sas98}), we arrive at the conclusion
that $U'={\mathcal W}_{\vec{\alpha},\vec{\beta}}U{\mathcal W}_{\vec{\alpha},\vec{\beta}}$.
Corollary \ref{esse50} is proved. \rule{2mm}{2mm}

\begin{corollary}\label{es60}
A  $J_{\vec{\alpha}}$-self-adjoint operator $A\in\Sigma_{J_{\vec{\alpha}}}$  commutes with ${J}_{\vec{\beta}}$, where
$\vec{\alpha}\cdot\vec{\beta}=0$ if and only if the corresponding unitary operator $U$ in (\ref{sas98})
satisfies the relation
\begin{equation}\label{es51}
 {\mathcal J}_{\vec{\beta}}U=U^{-1}{\mathcal J}_{\vec{\beta}}.
\end{equation}
\end{corollary}
\emph{Proof.}
Assume that $\vec{\beta}\in{\mathbb{S}^2}$ and  $\vec{\alpha}\cdot\vec{\beta}=0$.
Then $J_{\vec{\beta}}J_{\vec{\alpha}}=-J_{\vec{\alpha}}J_{\vec{\beta}}$ due to Lemma \ref{ll2}.
Since $J_{\vec{\beta}}S^*=S^*J_{\vec{\beta}}$, the commutation relation $AJ_{\vec{\beta}}=J_{\vec{\beta}}A$
is equivalent to the condition
\begin{equation}\label{esse70}
\forall{f}\in\mathcal{D}(A)\ {\Rightarrow} \ J_{\vec{\beta}}f\in\mathcal{D}(A).
\end{equation}

Let $f\in\mathcal{D}(S^*)$. Recalling that ${\mathcal {J}_{\vec{\beta}}}\Gamma_j=\Gamma_j{{J}_{\vec{\beta}}}$, we obtain
$$
\begin{array}{c}
({\mathcal J}_{\vec{\alpha}}\Gamma_1+i\Gamma_0)J_{\vec{\beta}}f=-{\mathcal {J}_{\vec{\beta}}}({\mathcal J}_{\vec{\alpha}}\Gamma_1-i\Gamma_0)f, \vspace{2mm} \\
({\mathcal J}_{\vec{\alpha}}\Gamma_1-i\Gamma_0)J_{\vec{\beta}}f=-{\mathcal {J}_{\vec{\beta}}}({\mathcal J}_{\vec{\alpha}}\Gamma_1+i\Gamma_0)f.
\end{array}
$$
Combining the last two relations with (\ref{sas98}), we conclude that (\ref{esse70}) is equivalent
to the identity
${\mathcal J}_{\vec{\beta}}U^{-1}{\mathcal J}_{\vec{\beta}}=U.$ Corollary \ref{es60} is proved. \rule{2mm}{2mm}

\begin{corollary}\label{es50}
A  $J_{\vec{\alpha}}$-self-adjoint operator $A\in\Sigma_{J_{\vec{\alpha}}}$  belongs to
the subset $\Upsilon$ (see (\ref{es9})) if and only if the corresponding unitary operator $U$ in (\ref{sas98})
satisfies the equality (\ref{es51}) \emph{for all} $\vec{\beta}\in{\mathbb{S}^2}$ such that $\vec{\alpha}\cdot\vec{\beta}=0$.
\end{corollary}
\emph{Proof.} Since $U$ satisfies (\ref{es51}) for all $\vec{\beta}\in{\mathbb{S}^2}$ such that $\vec{\alpha}\cdot\vec{\beta}=0$,
the operator $A\in\Sigma_{J_{\vec{\alpha}}}$ commutes with an arbitrary $J_{\vec{\beta}}$ such that $J_{\vec{\alpha}}J_{\vec{\beta}}=-J_{\vec{\beta}}J_{\vec{\alpha}}$ (due to Lemma \ref{ll2} and Corollary \ref{es60}).
In particular, the fundamental symmetry $J_{\vec{\gamma}}=iJ_{\vec{\alpha}}J_{\vec{\beta}}$ anti-commutes with $J_{\vec{\alpha}}$ and
hence, $J_{\vec{\gamma}}A=AJ_{\vec{\gamma}}$. On the other hand, since $A\in\Sigma_{J_{\vec{\alpha}}}$, we have $J_{\vec{\alpha}}A=A^*J_{\vec{\alpha}}$ and
$$
J_{\vec{\gamma}}A=iJ_{\vec{\alpha}}J_{\vec{\beta}}A=iJ_{\vec{\alpha}}AJ_{\vec{\beta}}=A^*iJ_{\vec{\alpha}}J_{\vec{\beta}}=A^*J_{\vec{\gamma}}.
$$
Thus $AJ_{\vec{\gamma}}=A^*J_{\vec{\gamma}}$ and hence, $A=A^*$. This means that the self-adjoint extension $A\supset{S}$ commutes with all fundamental symmetries from the Clifford algebra ${\mathcal C}l_2(J,R)$. Therefore, $A\in\Upsilon.$ The inverse statement follows from
Corollary \ref{es60}. \rule{2mm}{2mm}

\begin{corollary}\label{esse80}
Let $A\in\Sigma_{J_{\vec{\alpha}}}$ be defined by (\ref{sas98}) with $U={\mathcal J}_{\vec{\gamma}}$,
where $\vec{\gamma}\in{\mathbb{S}^2}$ is an arbitrary vector such that $\vec{\alpha}\cdot\vec{\gamma}=0$.
Then $\sigma(A)=\mathbb{C}$, i.e., $A$ has  empty resolvent set.
\end{corollary}
\emph{Proof.}
 Taking into account that ${\mathcal J}_{\vec{\alpha}}{\mathcal J}_{\vec{\gamma}}=-{\mathcal J}_{\vec{\gamma}}{\mathcal J}_{\vec{\alpha}}$ (since $\vec{\alpha}\cdot\vec{\gamma}=0$), we rewrite the definition (\ref{sas98}) of $A$:
\begin{equation}\label{sas98b}
A=S^*\upharpoonright\{f\in\mathcal{D}(S^*) \ : \ \mathcal{J}_{\vec{\alpha}}({\mathcal J}_{\vec{\gamma}}+I)\Gamma_1f=i({\mathcal J}_{\vec{\gamma}}+I)\Gamma_0f \}.
\end{equation}

Since relation (\ref{es51}) holds when $U={\mathcal J}_{\vec{\gamma}}$ and $\vec{\beta}=\vec{\gamma}$,
the operator $A\in\Sigma_{J_{\vec{\alpha}}}$ commutes with $J_{\vec{\gamma}}$ (Corollary \ref{es60}).
Therefore (cf. (\ref{esse5})),
\begin{equation}\label{esse5b}
A=\left(\begin{array}{cc}
A_{+} & 0 \\
0 & A_{-} \end{array}\right), \qquad A_{+}=A\upharpoonright_{\mathfrak{H}^{\gamma}_{+}}, \quad A_{-}=A\upharpoonright_{\mathfrak{H}_{-}^{\gamma}}
\end{equation}
with respect to the decomposition (\ref{es90}). Here
$S_{\gamma+}{\subseteq}A_{+}{\subseteq}S_{\gamma+}^*$ and $S_{\gamma-}{\subseteq}A_{-}{\subseteq}S_{\gamma-}^*$, where
$S_{\gamma\pm}=S\upharpoonright_{\mathfrak{H}_{\pm}^{\gamma}}$.

Denote $\mathcal{H}_{+}^{\gamma}=\frac{1}{2}(I+\mathcal{J}_{\vec{\gamma}})\mathcal{H}$ and
$\mathcal{H}_{-}^{\gamma}=\frac{1}{2}(I-\mathcal{J}_{\vec{\gamma}})\mathcal{H}$.  Then
\begin{equation}\label{usa14}
\mathcal{H}=\mathcal{H}_{+}^{\gamma}\oplus\mathcal{H}_{-}^{\gamma}
\end{equation}
and
$(\mathcal{H}_{\pm}^{\gamma}, \Gamma_0, \Gamma_1)$ are boundary triplets of operators $S^*_{\gamma\pm}$ (due to (\ref{es45}) and Lemma \ref{es16}).

Let $f\in\mathcal{D}(S^*_{\gamma+})$. Then $\Gamma_jf\in\mathcal{H}_{+}^{\gamma}$, $j=0,1$ and the identity in (\ref{sas98b}) takes the form
\begin{equation}\label{esse91}
\mathcal{J}_{\vec{\alpha}}\Gamma_1f=i\Gamma_0f.
\end{equation}
Since ${\mathcal J}_{\vec{\alpha}}{\mathcal J}_{\vec{\gamma}}=-{\mathcal J}_{\vec{\gamma}}{\mathcal J}_{\vec{\alpha}}$, the operator
${\mathcal J}_{\vec{\alpha}}$ maps $\mathcal{H}_{+}^{\gamma}$ onto $\mathcal{H}_{-}^{\gamma}$. Thus, (\ref{esse91}) may only hold in the case where $\Gamma_0f=\Gamma_1f=0$. Therefore, the operator $A_+$ in (\ref{esse5b}) coincides with $S_{+}^{\gamma}$.

Assume now $f\in\mathcal{D}(S^*_{\gamma-})$. Then $\Gamma_jf\in\mathcal{H}_{-}^{\gamma}$, $j=0,1$ and the identity in (\ref{sas98b})
vanishes (i.e., $0=0$). This means that $A_-=S^*_{\gamma-}$. Therefore, $A=A_\gamma$, where $A_\gamma$ is defined by (\ref{ew1}) and
$\sigma(A_\gamma)=\mathbb{C}$. \rule{2mm}{2mm}

\subsection{The resolvent formula.}
Let $\gamma(\mu)=(\Gamma_0\upharpoonright_{\ker(S^*-\mu{I})})^{-1}$ be the $\gamma$-field corresponding to the boundary triplet
$({\mathcal H}, \Gamma_0, \Gamma_1)$ of $S^*$ with properties (\ref{ea6a}).
Since ${J}_{\vec{\alpha}}$ maps $\ker(S^*-\mu{I})$ onto $\ker(S^*-\mu{I})$, formula (\ref{es45}) implies
$$
\gamma(\mu){\mathcal J}_{\vec{\alpha}}={J}_{\vec{\alpha}}\gamma(\mu), \quad \forall\mu\in\mathbb{C}\setminus\mathbb{R}
$$
for an arbitrary fundamental symmetry ${J}_{\vec{\alpha}}\in{\mathcal C}l_2({J},{R})$.

Let $A_0=S^*\upharpoonright\ker\Gamma_0$. Then $A_0$ is a self-adjoint extension of $S$ (due to the general properties of boundary triplets \cite{Gor}). Moreover, it follows from (\ref{es9}) and Remark \ref{es21} that $A_0\in\Upsilon$.

\begin{proposition}\label{neww909}
Let $(\mathcal{H}, \Gamma_0, \Gamma_1)$ be a boundary triplet
of $S^*$ with properties (\ref{ea6a}) and let $A\in\Sigma_{J_{\vec{\alpha}}}$ be defined by (\ref{sas98}).
Assume that $A$ is disjoint with $A_0$ (i.e., $\mathcal{D}(A)\cap\mathcal{D}(A_0)=\mathcal{D}(S)$) and
$\mu\in\rho(A)\cap\rho(A_0)$, then
\begin{equation}\label{bbb1}
(A-\mu{I})^{-1}=(A_0-\mu{I})^{-1}-\gamma(\mu)[M(\mu)-T]^{-1}\gamma^*(\overline{\mu}),
\end{equation}
where $T=i{\mathcal J}_{\vec{\alpha}}(I+U)(I-U)^{-1}$ is a ${\mathcal J}_{\vec{\alpha}}$-self-adjoint operator
in the Krein space $({\mathcal H}, [\cdot,\cdot]_{{\mathcal J}_{\vec{\alpha}}})$.
\end{proposition}
\emph{Proof.}
Since $A$ and $A_0$ are disjoint, the unitary operator $U$ which corresponds to the operator $A\in\Sigma_{J_{\vec{\alpha}}}$ in (\ref{sas98}) satisfies the relation $\ker(I-U)=\{0\}$. This relation and (\ref{es45}) allows one to rewrite (\ref{sas98}) as follows:
\begin{equation}\label{as1bbb}
A=S^*\upharpoonright{\{f\in\mathcal{D}(S^*) \ | \
T\Gamma_{0}f=\Gamma_{1}f\}},
\end{equation}
where $T=i{\mathcal J}_{\vec{\alpha}}(I+U)(I-U)^{-1}$ is a ${\mathcal J}_{\vec{\alpha}}$-self-adjoint operator
in the Krein space $({\mathcal H}, [\cdot,\cdot]_{{\mathcal J}_{\vec{\alpha}}})$ (due to self-adjointness of $i(I+U)(I-U)^{-1}$).
Repeating the standard arguments (see, e.g., \cite[p.14]{DM}), we deduce (\ref{bbb1}) from (\ref{as1bbb}).
\rule{2mm}{2mm}
\begin{remark}
The condition of disjointness of $A$ and $A_0$ in Proposition \ref{neww909} is not essential and it is
assumed for simplifying the exposition. In particular, this allows one to avoid  operators $A$ with empty resolvent set (see Corollary \ref{esse80} and relation (\ref{sas98b})) for which the formula (\ref{bbb1}) has no sense.
In the case of an arbitrary $A\in\Sigma_{J_{\vec{\alpha}}}$ with non-empty resolvent set,
the formula (\ref{bbb1}) also remains true  if we interpret $T$ as a $\mathcal{J}_{\vec{\alpha}}$-self-adjoint relation in $\mathcal{H}$
(see \cite[Theorem 3.22]{HK} for a similar result and \cite{BGP} for the basic definitions of linear relations theory).
\end{remark}

\section{The case of deficiency indices $<2,2>$}
In what follows, the symmetric operator $S$ has deficiency indices $<2,2>$.

{\bf 1.} Let  $(\mathcal{H}, \Gamma_0, \Gamma_1)$ be a boundary triplet of $S^*$ with properties
(\ref{ea6a}) or, that is equivalent, with properties (\ref{es45}). Let us fix an arbitrary fundamental symmetry
${\mathcal J}_{\vec{\gamma}}\in{\mathcal Cl}_2({\mathcal J}, {\mathcal R})$ and consider the decomposition
$\mathcal{H}=\mathcal{H}_{+}^{\gamma}\oplus\mathcal{H}_{-}^{\gamma}$
constructed by ${\mathcal J}_{\vec{\gamma}}$ (see (\ref{usa14})).
Then the Weyl function $M(\cdot)$ associated with $(\mathcal{H}, \Gamma_0, \Gamma_1)$ can be rewritten as
$$
M(\cdot)=\left(\begin{array}{cc}
m_{++}(\cdot) & m_{+-}(\cdot) \\
m_{-+}(\cdot) & m_{--}(\cdot)
\end{array}\right), \qquad m_{xy}(\cdot) : \mathcal{H}_y^\gamma\to\mathcal{H}_x^\gamma, \quad x, y\in\{+, -\},
$$
where $m_{xy}(\cdot)$ are scalar functions
(since $\dim{\mathcal H}=2$ and $\dim{\mathcal H}_\pm^\gamma=1$).

According to Lemma \ref{es30}, $M(\cdot)$ commutes with every fundamental symmetry from
${\mathcal C}l_2({\mathcal J},{\mathcal R})$. In particular, $\sigma_{j}M(\cdot)=M(\cdot)\sigma_{j}$ $(j=1,3)$,
where
$\sigma_1=\left(\begin{array}{cc} 0  & 1 \\
1 & 0
\end{array}\right)$, \ $\sigma_3=\left(\begin{array}{cc} 1  & 0 \\
0 & -1
\end{array}\right)$ are Pauli matrices.
This is possible only in the case
$$
m_{+-}(\cdot)=m_{-+}(\cdot)=0, \qquad m_{++}(\cdot)=m_{--}(\cdot),
$$
i.e.,
\begin{equation}\label{sas94}
M(\cdot)=m(\cdot)E,
\end{equation}
where $m(\cdot)=m_{++}(\cdot)=m_{--}(\cdot)$ is a scalar function defined on $\mathbb{C}\setminus\mathbb{R}$ and  $E$ is the identity $2\times{2}$-matrix.

Recalling that $(\mathcal{H}_{+}^{\gamma}, \Gamma_0, \Gamma_1)$ is a boundary triplet of $S^*_{\gamma+}$ (see the proof of Corollary \ref{esse80}) and taking into account the definition (\ref{neww65}) of Weyl functions, we obtain
the following statement which initially was established in \cite[Subsection 4.2]{HK}:
 \begin{proposition}
Let $({\mathcal H}, \Gamma_0, \Gamma_1)$ be a boundary triplet of $S$ defined above.  Then the Weyl function
$M(\cdot)$ is defined by (\ref{sas94}), where $m(\cdot)$ is the Weyl function of $S_{\gamma+}$ associated with boundary triplet
$(\mathcal{H}_{+}^{\gamma}, \Gamma_0, \Gamma_1)$. The function $m(\cdot)$ does not depend on the choice of $\vec{\gamma}\in\mathbb{S}^2$.
\end{proposition}

{\bf 2.} Let  $\vec{\alpha}, \vec{\beta}\in{\mathbb{S}^2}$ be
linearly independent vectors. According to Corollary \ref{usa3} all operators
$A\in\Sigma_{J_{\vec{\alpha}}}\cap\Sigma_{J_{\vec{\beta}}}$ are described by the formula
(\ref{usa8}). This means that spectra of these operators are completely characterized by
the spectra of their counterparts $A_+$ in (\ref{usa8}).

The operator $A_+$ is supposed to be an intermediate extension of $S_{\gamma+}$.
Two different situations may occur: 1. $A_+=S_{\gamma+}$ or $A_+=S_{\gamma+}^*$; \
2. $A_+$ is a quasi-self-adjoint extension\footnote{this class includes self-adjoint extensions also} of $S$, i.e., $S_{\gamma+}{\subset}A_{+}{\subset}S_{\gamma+}^*$. In the first case, the operators
$A\in\Sigma_{J_{\vec{\alpha}}}\cap\Sigma_{J_{\vec{\beta}}}$ have empty resolvent set (Theorem \ref{um2});
in the second case, the spectral properties of $A_+$ (and hence, $A$) are well known (see, e.g., \cite[Theorem 1, Appendix I]{AkGl}).
Summing up, we arrive at the following conclusion:
\begin{proposition}
Let $S$ be a simple symmetric operator with deficiency indices $<2,2>$
and $A\in\Sigma_{J_{\vec{\alpha}}}\cap\Sigma_{J_{\vec{\beta}}}$. Then or $\sigma(A)=\mathbb{C}$ or
the spectrum of $A$ consists of the spectral kernel of $S$ and the set of eigenvalues which can have
only real accommodation points.
\end{proposition}

{\bf 3.} Denote by $\Xi_{{\vec{\alpha}}}$ the collection of all operators $A\in\Sigma_{J_{\vec{\alpha}}}$ with empty resolvent set:
$$
\Xi_{{\vec{\alpha}}}=\{A\in\Sigma_{J_{\vec{\alpha}}} \ : \ \rho(A)=\emptyset \ \}
$$
and by $\Xi_{\vec{\alpha},\vec{\beta}}$ the pair of two operators $A_{\gamma(\vec{\beta})}$ and $A_{\gamma(\vec{\beta})}^*$
with empty resolvent set which are defined by (\ref{ew1}) for a fixed ${\vec{\alpha}}$ and ${\vec{\beta}}$.

\begin{theorem}\label{usa20}
Assume that $S$ is a symmetric operator with deficiency indices $<2,2>$ and its Weyl function
(associated with an arbitrary boundary triplet) differs from constant on $\mathbb{C}\setminus\mathbb{R}$.
Then
\begin{equation}\label{usa21}
\Xi_{{\vec{\alpha}}}=\bigcup_{\forall\vec{\beta}\in\mathbb{S}^2, \ \vec{\alpha}\cdot\vec{\beta}=0}\Xi_{\vec{\alpha},\vec{\beta}}.
\end{equation}
\end{theorem}
\emph{Proof.} By Theorem \ref{um2}, $\Xi_{{\vec{\alpha}}}\supset\Xi_{\vec{\alpha},\vec{\beta}}$ for all $\vec{\beta}\in\mathbb{S}^2$
such that $\vec{\alpha}\cdot\vec{\beta}=0$. Therefore, $\Xi_{{\vec{\alpha}}}\supset\bigcup\Xi_{\vec{\alpha},\vec{\beta}}$.

In the case of deficiency indices $<2,2>$ of $S$, the set $\Xi_{{\vec{\alpha}}}$ of all $J_{\vec{\alpha}}$-self-adjoint extensions
with empty resolvent set is described in \cite{KT}.
We briefly outline the principal results.

Denote by $
\mathfrak{N}_\mu=\ker(S^*-\mu{I})$, $\mu\in\mathbb{C}\setminus\mathbb{R}$,
the defect subspaces of $S$  and consider
the Hilbert space $\mathfrak{M}=\mathfrak{N}_{i}\dot{+}\mathfrak{N}_{-i}$ with
the inner product
$$
(x,y)_{\mathfrak{M}}=2[(x_i,y_i)+(x_{-i},y_{-i})],
$$
where $x=x_i+x_{-i}$ and $y=y_i+y_{-i}$ with  $x_{i}, y_{i}\in\mathfrak{N}_i$,
$x_{-i}, y_{-i}\in\mathfrak{N}_{-i}$.

The operator $Z$ that acts as identity
operator $I$ on $\mathfrak{N}_{i}$ and minus identity operator $-I$ on
$\mathfrak{N}_{-i}$ is an example of fundamental
symmetry in $\mathfrak{M}$. Other examples can be constructed due to the fact that $S$ commutes with $J_{\vec{\beta}}$
for all $\vec{\beta}\in\mathbb{S}^2$. This means that the
subspaces $\mathfrak{N}_{\pm{i}}$ reduce $J_{\vec{\beta}}$ and the restriction
$J_{\vec{\beta}}\upharpoonright\mathfrak{M}$ gives rise to a fundamental
symmetry in the Hilbert space $\mathfrak{M}$. Moreover, according to
the properties of $Z$ mentioned above, $J_{\vec{\beta}}Z=ZJ_{\vec{\beta}}$ and $J_{\vec{\beta}}Z$ is
a fundamental symmetry in $\mathfrak{M}$. Therefore, the
sesquilinear form
$$
[x,y]_{J_{\vec{\beta}}Z}=(J_{\vec{\beta}}Zx,y)_{\mathfrak{M}}=2[(J_{\vec{\beta}}x_i,y_i)-(J_{\vec{\beta}}x_{-i},y_{-i})]
$$
defines an indefinite metric on $\mathfrak{M}$.

According to the von-Neumann formulas, any closed intermediate extension
$A$ of $S$ (i.e.,
$S\subseteq{A}\subseteq{S}^*$) is
uniquely determined by the choice of a subspace $M\subset\mathfrak{M}$:
\begin{equation}\label{e55}
A=S^*\upharpoonright_{\mathcal{D}(A)}, \qquad \mathcal{D}(A)=\mathcal{D}(S)\dot{+}M,
\end{equation}

In particular, $J_{\vec{\beta}}$-self-adjoint extensions $A$ of $S$ correspond to
\emph{hypermaximal neutral subspaces $M$ with respect to  $[\cdot, \cdot]_{{J_{\vec{\beta}}Z}}$}.
This means that $A\in\Sigma_{J_{\vec{\alpha}}}\cap\Sigma_{J_{\vec{\beta}}}$ if and only if the corresponding subspace $M$ in (\ref{e55})
is \emph{simultaneously hypermaximal neutral with respect two different indefinite metrics $[\cdot, \cdot]_{{J_{\vec{\alpha}}Z}}$ and
$[\cdot, \cdot]_{{J_{\vec{\beta}}Z}}$}.

Without loss of generality we assume that $J_{\vec{\alpha}}$ coincides with $J$ in (\ref{es100b}), i.e.,
${\vec{\alpha}}=(1,0,0)$. Then fundamental symmetries $J_{\vec{\beta}}$ which anti-commute with $J$ have the form
\begin{equation}\label{usa22}
J_{\vec{\beta}}=\beta_2R+\beta_3iJR, \qquad \beta_2^2+\beta_3^2=1.
\end{equation}

To specify $M$ we consider an orthonormal basis $\{e_{++}, e_{+-}, e_{-+}, e_{--}\}$ of $\mathfrak{M}$
which satisfies the relations
\begin{equation}\label{ss1}
\begin{array}{c}
Ze_{++}=e_{++}, \ \ Ze_{+-}=e_{+-}, \ \ Ze_{-+}=-e_{-+}, \ \ Ze_{--}=-e_{--} \\[2mm]
Je_{++}=e_{++}, \ \ Je_{+-}=-e_{+-}, \ \ Je_{-+}=e_{-+}, \ \ Je_{--}=-e_{--} \\[2mm]
Re_{++}=e_{+-}, \ \ Re_{+-}=e_{++}, \ \  Re_{--}=e_{-+}, \ \
Re_{-+}=e_{--}.
\end{array}
\end{equation}

The existence of this basis was established in \cite{AKG} and it was used in \cite[Corollary 3.2]{KT}
to describe the collection of all $M$ in (\ref{e55}) which correspond to $J$-self-adjoint extensions of $S$ with
empty resolvent set. Such a description depends on properties of Weyl function of $S$. In particular, if
the Weyl function differs from the constant for a fixed boundary triplet, then this property remains
true for Weyl functions associated with an arbitrary boundary triplet of $S$.
Then, using relations (2.7)-(2.9) in \cite{KT}, we deduce that the Straus
characteristic function of $S$  (see \cite{SH}) differs from the zero-function on  $\mathbb{C}\setminus\mathbb{R}$.
In this case, Corollary 3.2 says that a $J$-self-adjoint extension  $A$ has empty resolvent set if and only if
the corresponding subspace $M$ coincides with linear
span $M={\textsf{span}}\{d_1, d_2\}$, where $d_1=e_{++}+e^{i\gamma}e_{+-}$, $d_2=e_{--}+e^{-i\gamma}e_{-+}$, and $\gamma\in[0,2\pi)$ is an arbitrary parameter.

The operator $A$ will belong to $\Sigma_{J_{\vec{\beta}}}$ if and only if the subspace $M={\textsf{span}}\{d_1, d_2\}$
turns out to be  hypermaximal neutral with respect to $[\cdot, \cdot]_{{J_{\vec{\beta}}Z}}$. Since $\dim{M}=2$
and $\dim\mathfrak{M}=4$, it suffices to check the neutrality of $M$. The last condition is equivalent to the relations
$$
[d_1, d_2]_{{J_{\vec{\beta}}Z}}=0, \qquad [d_1, d_1]_{{J_{\vec{\beta}}Z}}=0, \qquad [d_2, d_2]_{{J_{\vec{\beta}}Z}}=0.
$$

Using (\ref{usa22}), (\ref{ss1}), and remembering the orthogonality of $e_{\pm,\pm}$ in $\mathfrak{M}$, we
establish that $[d_1, d_2]_{{J_{\vec{\beta}}Z}}=0$ for all $\gamma\in[0,2\pi)$. The next two conditions are
transformed to the linear equation
\begin{equation}\label{usa25}
(\cos\gamma)\beta_2-(\sin\gamma)\beta_3=0,
\end{equation}
which has the nontrivial solution $\beta_2=\sin\gamma$,  $\beta_3=\cos\gamma$ for any $\gamma\in[0,2\pi)$. This means that
an arbitrary $J$-self-adjoint extension $A$ with empty resolvent set is also
a $J_{\vec{\beta}}$-self-adjoint operator under choosing $\beta_2$ and $\beta_3$
in (\ref{usa22}) as solutions of (\ref{usa25}).
Theorem \ref{usa20} is proved.
\rule{2mm}{2mm}

\begin{corollary}\label{usa26}
Let $S$ be a symmetric operator with deficiency indices $<2,2>$ and let $(\mathcal{H}, \Gamma_0, \Gamma_1)$ be a boundary triplet
of $S^*$ with properties (\ref{ea6a}). If the Weyl function of $S$ differs from constant on $\mathbb{C}\setminus\mathbb{R}$, then
the set $\Xi_{{\vec{\alpha}}}$  is described by
(\ref{sas98}) where $U$ runs the set of all fundamental symmetries ${\mathcal J}_{\vec{\beta}}\in{\mathcal C}l_2(\mathcal{J},\mathcal{R})$
such that ${\vec{\alpha}}\cdot{\vec{\beta}}=0$.
\end{corollary}
\emph{Proof.} It follows from Corollary \ref{esse80} and Theorem \ref{usa20}. \rule{2mm}{2mm}

Theorem \ref{usa20} and Corollary \ref{usa26} are not true when the Weyl function of $S$ is a constant.
In that case, the set $\Xi_{{\vec{\alpha}}}$ of $J_{\vec{\alpha}}$-self-adjoint extensions increases considerably and
$\Xi_{{\vec{\alpha}}}\supset\bigcup\Xi_{\vec{\alpha},\vec{\beta}}$.

\begin{corollary}\label{usa70}
Let $S$ be a simple symmetric operator with deficiency indices $<2,2>$. Then the following statements are equivalent:
\begin{itemize}
  \item[(i)] the strict inclusion
$$
\Xi_{{\vec{\alpha}}}\supset\bigcup_{\forall\vec{\beta}\in\mathbb{S}^2, \ \vec{\alpha}\cdot\vec{\beta}=0}\Xi_{\vec{\alpha},\vec{\beta}}
$$
holds;
\item[(ii)] the Weyl function $M(\cdot)$ of $S$ is a constant on $\mathbb{C}\setminus\mathbb{R}$;
\item[(iii)] $S$ is unitarily equivalent to the symmetric operator in $L_2(\mathbb{R}, \mathbb{C}^2)$:
\begin{equation}\label{usa72}
S'=i\frac{d}{dx}, \qquad  \mathcal{D}(S')=\{u\in{W_2^1(\mathbb{R}, \mathbb{C}^2)} \ : \ u(0)=0\}.
\end{equation}
\end{itemize}
 \end{corollary}
\emph{Proof.} Assume that the Weyl function $M(\cdot)$ of $S$ is a constant. By (\ref{sas94}), the Weyl function $m(\cdot)$ of $S_{\gamma+}=S\upharpoonright_{\mathfrak{H}_+^\gamma}$ is also constant.
This means that the Straus characteristic function of the simple symmetric operator $S_{\gamma+}$
with deficiency indices $<1, 1>$ is zero on  $\mathbb{C}\setminus\mathbb{R}$ (see the proof of Theorem \ref{usa20}).
Therefore, $S_{\gamma+}$ is unitarily equivalent to the symmetric operator $S'_+=i\frac{d}{dx}, \ \mathcal{D}(S'_+)=\{u\in{W_2^1(\mathbb{R})} \ : \ u(0)=0\}$ in $L_2(\mathbb{R})$ \cite[Subsection 3.4]{KU}.

Recalling the decomposition (\ref{es1b}) of $S$, where the simple symmetric operator $S_{\gamma-}=S\upharpoonright_{\mathfrak{H}_-^\gamma}$ also has  deficiency indices $<1, 1>$  and  zero characteristic function, we conclude that $S$ is unitarily equivalent to the symmetric operator $S'$
defined by (\ref{usa72}). This establishes the equivalence of (ii) and (iii).

Assume again that the Weyl function of $S$ is a constant. Then the Straus characteristic function of $S$ is zero. In that case,
Corollary 3.2 in \cite{KT} yields that $A\in\Xi_{{\vec{\alpha}}}$ if and only if
the corresponding subspace $M$ in (\ref{e55}) coincides with linear
span
\begin{equation}\label{usa80}
M={\textsf{span}}\{d_1, d_2\}, \quad d_1=e_{++}+e^{i(\phi+\gamma)}e_{+-}, \quad d_2=e_{--}+e^{i(\phi-\gamma)}e_{-+},
\end{equation}
where  $\phi, \gamma\in[0,2\pi)$ are \emph{two arbitrary parameters}. Thus the set $\Xi_{{\vec{\alpha}}}$ is described by
two independent parameters $\phi$ and  $\gamma$.

Due to the proof of Theorem \ref{usa20}, the operator $A\in\Xi_{{\vec{\alpha}}}$ belongs to the subset $\Xi_{\vec{\alpha},\vec{\beta}}$ if and only if the subspace $M$ in (\ref{usa80}) is neutral with respect to $[\cdot, \cdot]_{{J_{\vec{\beta}}Z}}$.  Repeating the argumentation above, we conclude that the neutrality of $M$ is  equivalent to the existence of nontrivial solution $\beta_2, \beta_3$
of the system (cf. (\ref{usa25}))
\begin{equation}\label{usa81}
\left\{\begin{array}{l}
\cos(\phi+\gamma)\beta_2-\sin(\phi+\gamma)\beta_3=0 \vspace{1mm} \\
\cos(\phi-\gamma)\beta_2+\sin(\phi-\gamma)\beta_3=0
\end{array}\right.
\end{equation}

The determinant of (\ref{usa81}) is $\sin2\phi$. Therefore, there are no nontrivial solutions for $\phi\not\in\{0, \frac{\pi}{2}, \pi, \frac{3\pi}{2}\}$. This means the existence of operators $A\in\Xi_{{\vec{\alpha}}}$ which, simultaneously, do not belong to $\bigcup\Xi_{\vec{\alpha},\vec{\beta}}$.
Thus, we establish the equivalence of (ii) and (i). Corollary \ref{usa70} is proved
 \rule{2mm}{2mm}

{\bf 4.} Consider the one-dimensional Schr\"{o}dinger differential expression
\begin{equation}\label{usa51}
l(\phi)(x)=-\phi''(x)+q(x)\phi(x), \qquad x\in\mathbb{R},
\end{equation}
where $q$ is an \emph{even} real-valued measurable function that has a non-integrable singularity at zero and is integrable on every finite subinterval of $\mathbb{R}\setminus\{0\}$.

Assume in what follows that the potential $q(x)$ is in the limit point case at $x\to\pm\infty$ and is in the limit-circle case
at $x=0$. Denote by $\mathcal{D}$ the set of all functions $\phi(x)\in{L_2(\mathbb{R})}$ such that $\phi$ and $\phi'$ are absolutely continuous on every finite subinterval of $\mathbb{R}\setminus\{0\}$ and $l(\phi)\in{L_2(\mathbb{R})}$. On $\mathcal{D}$ we define the operator
$L$ as follows:
$$
L\phi=l(\phi), \qquad \forall\phi\in\mathcal{D}.
$$

The operator $L$ commutes with the space parity operator  $\mathcal{P}\phi(x)=\phi(-x)$ and with the operator of multiplication by
 $(\textsf{sgn}\ x)I$. These operators are anti-commuting fundamental symmetries in $L_2(\mathbb{R})$.
Therefore, $L$ commutes with elements of the Clifford algebra $\mathcal{C}l_2(\mathcal{P}, (\textsf{sgn}\ x)I)$.
However, $L$ is not a symmetric operator.

Denote for brevity  $J_{\vec{\gamma}}=(\textsf{sgn}\ x)I$. Then, the decomposition
(\ref{es90}) takes the form $L_2(\mathbb{R})=L_2(\mathbb{R}_+){\oplus}L_2(\mathbb{R}_-)$ and with respect to it
\begin{equation}\label{usa30}
L=\left(\begin{array}{cc}
L_+ & 0 \\
0 & \mathcal{P}L_+\mathcal{P}
\end{array}\right), \qquad L_+=L\upharpoonright_{L_2(\mathbb{R}_+)}.
\end{equation}

The operator $L_+$ is the maximal operator for differential expression $l(\phi)$ considered on semi-axes $\mathbb{R}_+=(0,\infty)$.
Denote by $S_+$ the minimal operator generated $l(\phi)$ in $L_2(\mathbb{R}_+)$. The symmetric operator $S_+$ has
deficiency indices $<1,1>$.

Let $(\mathbb{C}, \Gamma_0^+, \Gamma_1^+)$ be an arbitrary boundary
triplet of $L_+=S_+^*$ in $L_2(\mathbb{R}_+)$. Then, the boundary triplet $(\mathbb{C}^2, \Gamma_0, \Gamma_1)$
determined by (\ref{esse41}) with $J_{\vec{\tau}}=\mathcal{P}$ and $N=\mathbb{C}$ is a boundary triplet of $L=S^*$
in the space $L_2(\mathbb{R})$. Here $S=\left(\begin{array}{cc}
S_+ & 0 \\
0 & \mathcal{P}S_+\mathcal{P}
\end{array}\right)$ is the symmetric operator in $L_2(\mathbb{R})$ (cf. (\ref{usa30})) with deficiency indices $<2,2>$.

Let $J_{\vec{\alpha}}$ be an arbitrary fundamental symmetry from ${\mathcal Cl}_2(\mathcal{P}, (\textsf{sgn}\ x)I)$.
By Theorem \ref{new45}, $J_{\vec{\alpha}}$-self-adjoint extensions $A\in\Sigma_{J_{\vec{\alpha}}}$ of $S$  are defined as the restrictions
of $L$:
\begin{equation}\label{sas98ba}
A=L\upharpoonright\{f\in\mathcal{D} \ : \ U(\mathcal{J_{\vec{\alpha}}}\Gamma_1+i\Gamma_0)f=(\mathcal{J_{\vec{\alpha}}}\Gamma_1-i\Gamma_0)f \},
\end{equation}
where $U$ runs the set of $2\times{2}$-unitary matrices.
The operators $A$ can be interpreted as $J_{\vec{\alpha}}$-self-adjoint operator realizations of differential
expression (\ref{usa51}) in $L_2(\mathbb{R})$.

Since the sets $\Sigma_{J_{\vec{\alpha}}}$ are unitarily equivalent
for different ${\vec{\alpha}}\in\mathbb{S}^2$ (Theorem \ref{es12a}) one can set  $J_{\vec{\alpha}}=\mathcal{P}$ for definiteness.
\begin{proposition}\label{usa90}
The collection of all $\mathcal{P}$-self-adjoint extensions $A\in\Sigma_{\mathcal P}$ with empty resolvent set
coincides with the restrictions of $L$ onto the sets of functions $f\in\mathcal{D}$  satisfying the condition
$$
\left(\begin{array}{cc}
i\sin\theta & 1-\cos\theta \\
1+\cos\theta & -i\sin\theta
\end{array}\right)\Gamma_1f=i\left(\begin{array}{cc}
1+\cos\theta & -i\sin\theta  \\
i\sin\theta & 1-\cos\theta
\end{array}\right)\Gamma_0f, \quad \forall\theta\in[0, 2\pi).
$$
\end{proposition}
\emph{Proof.} Since $J_{\vec{\tau}}={\mathcal P}$ and $J_{\vec{\gamma}}=(\textsf{sgn}\ x)I$, relations (\ref{ea6a}), (\ref{usa31}) mean that
$\sigma_1\Gamma_j=\Gamma_j{\mathcal P}$ and
$\sigma_3\Gamma_j=\Gamma_j(\textsf{sgn}\ x)I$  $(j=0,1)$,
where $\sigma_1$ and $\sigma_3$ are Pauli matrices. Therefore, the `image' of
the Clifford algebra $\mathcal{C}l_2(\mathcal{P}, (\textsf{sgn}\ x)I)$ coincides with
$\mathcal{C}l_2(\sigma_1, \sigma_3)$ in $\mathbb{C}^2$ (see Remark \ref{es21}).

The fundamental symmetries ${\mathcal J}_{\vec{\beta}}\in{\mathcal C}l_2(\sigma_1, \sigma_3)$ anti-commuting with $\sigma_1$
have the form $
{\mathcal J}_{\vec{\beta}}=\beta_2\sigma_2+\beta_3\sigma_3$, \ $\beta_2^2+\beta_3^2=1,$
where  $\sigma_2=i\sigma_1\sigma_3$. Hence,
\begin{equation}\label{usa61}
{\mathcal J}_{\vec{\beta}}=\left(\begin{array}{cc}
\beta_3 & -i\beta_2 \\
i\beta_2 & -\beta_3
\end{array}\right)=\left(\begin{array}{cc}
\cos\theta & -i\sin\theta \\
i\sin\theta & -\cos\theta
\end{array}\right),  \quad \theta\in[0,2\pi).
\end{equation}
Here, we set $\beta_3=\cos\theta$ and $\beta_2=\sin\theta$ \ (since $\beta_2^2+\beta_3^2=1$).
Applying Corollary \ref{usa26} and rewriting $(\ref{sas98})$ in the form (\ref{sas98b}) with ${\mathcal J}_{\vec{\alpha}}=\sigma_1,$
${\mathcal J}_{\vec{\gamma}}={\mathcal J}_{\vec{\beta}}$ (here ${\mathcal J}_{\vec{\beta}}$ is determined by (\ref{usa61})), we complete
the proof of Proposition \ref{usa90}. \rule{2mm}{2mm}

\begin{remark}
To apply  Proposition \ref{usa90} for concrete potentials $q(x)$ in (\ref{usa61}) one needs only to construct
a boundary triplet $(\mathbb{C}^2, \Gamma_0, \Gamma_1)$ of $L$  with the help of  a boundary triplet
$(\mathbb{C}, \Gamma_0^+, \Gamma_1^+)$ of the differential expression (\ref{usa51}) on semi-axis $\mathbb{R}_+$ (see (\ref{esse41})).
To do that one can use \cite{KOCH}, where simple explicit formulas for operators $\Gamma_j^+$ constructed in terms of asymptotic behavior of $q(x)$ as $x\to{0}$ were obtained for great number of singular potentials.
\end{remark}

\noindent \textbf{Acknowledgements.}
The first author (S.K.) expresses his gratitude to the DFG
(project AL 214/33-1) and JRP IZ73Z0 (28135) of SCOPES 2009-2012 for the support.

\end{document}